\documentclass[12pt]{amsart}
\usepackage[top=1in, bottom=1in, left=1in, right=1in]{geometry}

\usepackage{amssymb,amsmath,amsthm}
\usepackage{mathrsfs} 
\usepackage{bbm} 
\usepackage{enumitem}
\usepackage[hyphens]{url}
\usepackage{times}

\usepackage{graphicx}
\usepackage{color}
\usepackage[colorlinks=true]{hyperref}

\numberwithin{equation}{section}

\newtheorem{thm}{Theorem}[section]
\newtheorem{lem}[thm]{Lemma}

\newtheorem{cor}[thm]{Corollary}

\theoremstyle{definition}

\theoremstyle{remark}
\newtheorem{remark}{Remark}[section]
\newtheorem*{remark*}{Remark}

\newtheoremstyle{claim} 
    {1em}                        
    {1em}                        
    {}                           
    {}                           
    {\bfseries}                  
    {.}                          
    {.5em}                       
    {}  
    
\theoremstyle{claim}    
\newtheorem{claim}{Claim}

\newcommand{\N}{\mathbb{N}}

\newcommand{\CA}{\mathcal{A}}

\newcommand{\CP}{\mathcal{P}}

\newcommand{\CS}{\mathcal{S}}

\newcommand{\abs}[1]{\left\lvert #1 \right\rvert}
\newcommand{\set}[1]{\left\{ #1 \right\}}

\newcommand{\bs}\boldsymbol{}

\renewcommand{\phi}{\varphi}

\renewcommand{\mod}[1]{\,({\rm mod}\,#1)}

\definecolor{red}{rgb}{1,0,0}

\definecolor{orange}{rgb}{0.7,0.3,0}

\definecolor{blue}{rgb}{.2,.6,.75}

\definecolor{green}{rgb}{.4,.7,.4}

\begin{document}

\title[Multiplicative Functions on Shifted Primes]{Multiplicative Functions on Shifted Primes}

\author{Stelios Sachpazis}
\address{D\'epartement de math\'ematiques et de statistique\\
Universit\'e de Montr\'eal\\
CP 6128 succ. Centre-Ville\\
Montr\'eal, QC H3C 3J7\\
Canada}
\email{stelios.sachpazis@umontreal.ca}


\date{\today}

\begin{abstract}
Let $f$ be a positive multiplicative function and let $k\geq 2$ be an integer. We prove that if the prime values $f(p)$ converge to $1$ sufficiently slowly as $p\rightarrow +\infty$, in the sense that $\sum_{p}|f(p)-1|=\infty$, there exists a real number $c>0$ such that the $k$-tuples $(f(p+1),\ldots,f(p+k))$ are dense in the hypercube $[0,c]^k$ or in $[c,+\infty)^k$. In particular, the values $f(p+1),\ldots,f(p+k)$ can be put in any increasing order infinitely often. Our work generalises previous results of De Koninck and Luca.
\end{abstract}

\maketitle


\section{Introduction}

\vspace{1mm}
 
The understanding of the local behavior of arithmetic functions has been the subject of research of many mathematicians. Part of this research involves the study of the values of an arithmetic function on consecutive integers. For example, if we denote the divisor function by $\tau$, in $1952$ Erd\"{o}s and Mirsky \cite{ermir} asked whether the equation $\tau(n)=\tau(n+1)$ admits infinitely many solutions in the set of natural numbers, a question that can be considered as a close relative of the twin prime conjecture. It remained open for about thirty years until Heath-Brown \cite{hb} answered it affirmatively in $1984$ by showing that 

$$\#\set{n\leq x:\tau(n)=\tau(n+1)}\gg x(\log x)^{-7}.$$

\vspace{3mm} 

The method of Heath-Brown also yielded that there exist infinitely many positive integer solutions to the equation $\Omega(n)=\Omega(n+1)$, where $\Omega$ counts the number of prime factors of $n$ with multiplicity. His method, however, was not successful in proving that the analogous equation $\omega(n)=\omega(n+1)$ has infinitely many positive integer solutions for the function $\omega$ which counts the number of distinct prime factors of $n$. It was Schlage-Puchta \cite{sp} who proved for the first time that even the equality $\omega(n)=\omega(n+1)$ holds for infinitely many natural numbers $n$. 

In $2011$, Goldston, Graham, Pintz and Yildirim \cite{mg} made a significant breakthrough. They did not just prove that there are infinitely many integer solutions of those equations, but they also showed that the value of the relevant arithmetic function $\tau,\ \Omega$ or $\omega$ can be specified. For instance, they proved that there are infinitely many integers $n$ such that $\omega(n)=\omega(n+1)=3.$

Arithmetic functions, such as $\tau$, are very sensitive on the exact number of prime factors of their input. In particular, they are highly sensitive on the large prime factors of their input. On the other hand, if we consider functions that are less sensitive to large prime factors, we can say more. Indeed, recently, De Koninck and Luca \cite{paper} proved that for any fixed integer $k\geq2$ we have that  

$$\frac{\max{\set{\phi(n+1),\ldots,\phi(n+k)}}}{\min{\set{\phi(n+1),\ldots,\phi(n+k)}}}$$

\vspace{3mm}

\noindent
is arbitrarily close to 1 infinitely often. Here $\phi$ stands for the Euler totient function. They also established the same result for the additive functions $\omega$, $\Omega$, the sum of divisors function $\sigma_1$ and the kernel function $\gamma$, which is the multiplicative function defined by the relation $\gamma(p^m)=p$ for any prime power $p^m$.

They had also proved \cite[Problem 8.6]{dl} that if $\set{i_1,\ldots,i_k}$ is a permutation of $\set{1,\ldots,k}$, the inequalities

\begin{equation}\label{in}
\phi(n+i_1)<\ldots<\phi(n+i_k)
\end{equation}

\vspace{1mm}

\noindent
hold for an infinite set of natural numbers $n.$

The purpose of the present paper is to extend these results. We state below our main result and an immediate consequence of it. Both of them are in the same spirit as the work of De Koninck and Luca.
 
\vspace{3mm}

\begin{thm}\label{bthm}
Let $k$ be a positive integer and let $f$ be a positive multiplicative function such that $\sum_{p}|f(p)-1|=\infty$ and $f(p)\rightarrow 1$ when $p\rightarrow \infty$. There exists a positive real number $c$ such that the set of tuples $(f(p+1),...,f(p+k))$ is dense either in $[0,c]^k$ or in $[c,+\infty)^k.$
\end{thm}

\vspace{1mm}

\begin{cor}\label{thm2}
Let $k\geq 2$ be an integer and let f be a positive multiplicative function such that $\sum_{p}|f(p)-1|=\infty$ and $f(p)\rightarrow 1$ when $p\rightarrow \infty$. If $\set{i_1,\ldots,i_k}$ is a permutation of $\set{1,\ldots,k}$, the inequalities

\begin{equation*}
f(p+i_1)<\ldots<f(p+i_k)
\end{equation*}

\vspace{1mm}

\noindent	
hold for infinitely many primes $p$.
\end{cor}

\begin{proof}
In Theorem \ref{bthm}, we may choose sufficiently small or large (this depends on whether we have density in $[0,c]^k$ or $[c,+\infty)^k$) non-negative numbers $c_1,\ldots,c_k$ such that $c_{i_1}<\ldots<c_{i_k}$ and $f(p+i_j)\rightarrow c_{i_j}$ for all $j\in \{1,\ldots,k\}$ on a subsequence of the primes. Then, for $\epsilon>0$ such that $\epsilon<\min{\set{(c_{i_{j+1}}-c_{i_j})/2,\ j=1,\ldots,k-1}}$, we have that

\begin{equation*}
f(p+i_j)<c_{i_j}+\epsilon<c_{i_{j+1}}-\epsilon<f(p+i_{j+1})
\end{equation*}

\vspace{2mm}

\noindent
for infinitely many primes $p$.
\end{proof}

\vspace{1mm}

\begin{remark}
	In Theorem \ref{bthm} we assumed that the function $f$ is positive. However, we can have an analogous result for $f$ being zero only on a finite set of primes by replacing the shifts $i\in\{1,\ldots,k\}$ by the $k$ first multiples of the product of those primes. For example, if $f(2)=0$, we can modify the proofs below by replacing the shifts $i\in \{1,\ldots,k\}$ by the even shifts $2i$. Doing so implies that there exists a constant $c>0$ such that the set of tuples $(f(p+2),\ldots,f(p+2k))$ is dense either in $[0,c]^k$ or $[c,+\infty)^k$.
\end{remark}

\vspace{1mm}

\begin{remark}
    This theorem extends (1.1), proven by De Koninck and Luca. It is also worth noting that there are additive functions whose values on consecutive integers can be ordered. For example, De Koninck, Friedlander and Luca \cite{dfk} proved that the inequalities
	
	\begin{equation}\label{str}
	\omega(n+1)<\ldots<\omega(n+k)
	\end{equation}
	
	\vspace{2mm}
	
	\noindent
	hold infinitely often. They also proved that $\omega$ can be replaced by $\Omega$ in (\ref{str}). Furthermore, minor changes in the solution of Problem $7.26$ in \cite{dl} can lead to these results with $\set{i_1,\ldots,i_k}=\set{1,\ldots,k}.$  
\end{remark}

\vspace{1mm}

Theorem \ref{bthm} may be directly deduced from its more technical analogue, which is the following theorem that we prove in Section 2.

\vspace{2mm}

\begin{thm}\label{prop}
Let $k\geq 2$ be an integer and let f be a positive multiplicative function such that $f(p)\rightarrow 1$ as $p\rightarrow \infty$. Let also $K:=k(k+1)((k-1)!)^2$.
	
(a) If
	
\begin{equation}\label{cond1}
\sum_{f(p)>1}(f(p)-1)=\infty,
\end{equation}
	
\noindent
then the set of tuples $(f(p+1),...,f(p+k))$ is dense in $[f(K),\infty)\times [f(1),\infty)\times\cdots\times[f(k-1),\infty).$
	
(b) If
	
\begin{equation}\label{cond2}
\sum_{f(p)<1}(1-f(p))=\infty,
\end{equation}
	
\noindent
the set of tuples $(f(p+1),...,f(p+k))$ is dense in $[0,f(K)]\times [0,f(1)]\times\cdots\times[0,f(k-1)].$
\end{thm}

\vspace{1mm}

\subsection*{Notation and definitions}

If $n$ is a positive integer such that $n\geq 2,$ then $P^+(n)$ denotes the largest prime factor of $n$ and $P^-(n)$ denotes the smallest prime factor of $n$. We also define $P^+(1)=1$ and $P^-(1)=+\infty.$ 

Let $y>0$ and let $n$ be a natural number. The $y$-rough part of $n$ is defined to be equal to $\prod_{p^a\|n,\, p>y}p^a$ and its $y$-smooth part is given by $\prod_{p^a\| n,\,p\le y}p^a$.

\vspace{3mm}

\subsection*{Acknowledgements}

The author would like to thank his advisor Dimitris Koukoulopoulos for all the useful discussions on Theorem \ref{prop} and for suggesting a sieving argument that led to a simplification of Theorem \ref{bthm}. He would also like to thank Andrew Granville for sharing his ideas about the necessity of the conditions in the statement of the main theorem. Last, but not least, he thanks the Stavros Niarchos Scholarships Foundation for the generous financial support that provides to him.

\vspace{2mm}

\section{Auxiliary Results and Proof of Theorem \ref{prop}}

\vspace{1mm}

In this section we prove two preparatory lemmas and use them to establish Theorem \ref{prop}. We begin with the first lemma, which will be only needed for the proof of the second one and we then make use of the second lemma to prove Theorem \ref{prop}.

\vspace{1mm}

\begin{lem}\label{lem2}
Let $\set{b_n}_{n\in \N}$ be a sequence of non-negative terms with $\sum_{n\geq 1}b_n=\infty$ and $b_n\rightarrow 0$ as $n\rightarrow \infty$. If $\beta>0,$ there exists a subsequence $\set{b_{n_k}}_{k\in \N}$ of $\set{b_n}_{n\in \N}$ such that $\sum_{k\in \N}b_{n_k}=\beta$.
\end{lem}

\begin{proof}
Since $b_n\rightarrow 0$, we may define $n_1$ to be the smallest positive integer such that $b_n<\beta$ for all $n\geq n_1.$ Moreover, since $\sum_{n\geq n_1}b_n=\infty,$ we may define $n_2$ to be the largest integer such that $n_2\geq n_1$ and 
	
$$S_1:=\sum_{n=n_1}^{n_2}b_n<\beta.$$
	
\noindent
Now, we inductively define $n_{2k-1},n_{2k}$ and $S_k$ for $k\geq 2.$ Since $b_n\rightarrow 0$, we define $n_{2k-1}$ to be the smallest integer such that $n_{2k-1}\geq n_{2(k-1)}$ and $b_{n}<\min{\set{\beta/k,\beta-(S_1+\ldots+S_{k-1})}}$ for all $n\geq n_{2k-1}$. Then, $n_{2k}$ is defined as the largest integer such that $n_{2k}\geq n_{2k-1}$ and
	
$$S_k:=\sum_{n=n_{2k-1}}^{n_{2k}}b_n<\beta-(S_1+\ldots+S_{k-1}).$$

By definition it is clear that $\beta-(S_1+\ldots+S_{k-1})\leq S_k+b_{n_{2k}+1}<S_k+\beta/k$. Consequently, $\beta(1-1/k)<\sum_{i=1}^kS_i<\beta.$ Now the proof is complete, because these inequalities imply that $\sum_{i\geq1}S_i=\beta$.
\end{proof}

\vspace{3mm}

Now that the proof of Lemma \ref{lem2} is complete, we continue by proving the second lemma of the section.

\vspace{3mm}

\begin{lem}\label{lem1}
Let $m\in \N$ and let $f$ be a multiplicative function with $f(p)\rightarrow 1$ as $p\rightarrow \infty$. 

(a) If $\epsilon>0$, $\set{x_n}_{n\in \N}\subset [1,+\infty)$ and $\sum_{f(p)>1}(f(p)-1)=\infty$, there exists a sequence of pairwise coprime positive integers $\{a_n\}_{n\in \mathbb{N}}$ such that 

\begin{equation*}
|f(a_n)-x_n|<\epsilon\quad \text{and}\quad (a_n,m)=1\quad \text{for all}\quad n\in \mathbb{N}.
\end{equation*}

\vspace{2mm}

(b) If $\epsilon>0$, $\set{x_n}_{n\in \N}\subset [0,1]$ and $\sum_{f(p)<1}(1-f(p))=\infty$, there exists a sequence of pairwise coprime positive integers $\{a_n\}_{n\in \mathbb{N}}$ such that 

\begin{equation*}
|f(a_n)-x_n|<\epsilon\quad \text{and}\quad (a_n,m)=1\quad \text{for all}\quad n\in \mathbb{N}.
\end{equation*}

\end{lem}

\begin{proof}
(a) First we assume that (\ref{cond1}) holds. Since $f(p)\rightarrow1$, there exists a prime $p_0>P^+(m)$ such that $|f(p)-1|\leq1/2$ for all $p\geq p_0.$ Then, 
	
\begin{equation}\label{ineq}
|f(p)-1|\asymp |\log f(p)| \quad \text{for all} \quad p\geq p_0
\end{equation} 
	
\noindent
and so, the series
	
$$\sum_{p\geq p_0,\ f(p)>1}\!\!\!\!\log f(p)$$
	
\noindent
diverges. 

Its terms are positive and converge to $0$. Therefore, if $x_1>1,$ we use Lemma \ref{lem2} with $\beta=\log x_1$ and find a square-free positive integer $a_1$ coming from the terms $f(p)>1$ with $p\geq p_0$ whose sums of logarithms are in $(\log(x_1-\epsilon),\log(x_1+\epsilon))$. Then $P^-(a_1)>P^+(m)$ and $|f(a_1)-x_1|<\epsilon.$ If $x_1=1$, then we define $a_1$ differently. Since $f(p)\rightarrow1$, there exists a prime $p_1>P^+(m)$ such that $|f(p_1)-1|<\epsilon$ and in this case we simply take $a_1=p_1.$

We continue by inductively defining $a_n$ for $n\geq 2.$ If $a_1,\ldots,a_{n-1}$ are already determined and $x_n>1$, we apply Lemma \ref{lem2} to the divergent series

$$\sum_{\substack{f(p)>1,\\p>P^+(ma_1\cdots a_{n-1})}}\!\!\!\!\!\!\!\!\log f(p)$$
	
\noindent
with $\beta=\log x_n.$ Consequently, there exists, as before, a positive square-free integer $a_n$ such that $|f(a_n)-x_n|<\epsilon$ and $P^-(a_n)>P^+(ma_1\cdots a_{n-1})$. If $x_n=1,$ we simply take $a_n=p_n$, where $p_n$ is a prime such that $p_n>P^+(ma_1\cdots a_{n-1})$ and $|f(p_n)-1|<\epsilon.$ The existence of such a prime is guaranteed by the fact that $f(p)\rightarrow1.$ Therefore, the proof of Lemma \ref{lem1} is complete when (\ref{cond1}) holds. 

(b) For the part (b) where (\ref{cond2}) holds, we can adjust the argument above to the divergent series
	
$$\sum_{f(p)<1}\log \left(\frac{1}{f(p)}\right).$$
	
\noindent
Note that its terms are positive and converge to $0.$ Its divergence follows using (\ref{ineq}).
\end{proof}

\vspace{3mm}

We close the section with the proof of Theorem \ref{prop}, which will be based on an auxiliary estimate (inequality (2.5)) that we will prove in the next and final section of the paper.

\vspace{3mm}

\begin{proof}[Proof of Theorem \ref{prop}]
Let $\nu \in \N$ and let $x$ be a real number $\geq 1$. The proof remains the same under each of the conditions (\ref{cond1}) and (\ref{cond2}). For $K:=k(k+1)((k-1)!)^2$, if (\ref{cond1}) holds, we let $c_1\in [f(K),+\infty)$ and $c_i\in [f(i-1),+\infty)$ for $i\in \{2,\ldots,k\}$, whereas if (\ref{cond2}) holds, we let $c_1\in[0,f(K)]$ and $c_i\in[0,f(i-1)]$ for $i\in \{2,\ldots,k\}$. For the proof it suffices to show that the cardinality of the set 
	
$$S(x,\nu):=\set{p\leq x: c_i-1/\nu<f(p+i)<c_i+1/\nu,\ i=1,\ldots,k}$$

\vspace{3mm}

\noindent
tends to infinity as $x\rightarrow \infty.$ We show this by constructing a subset of $S(x,\nu)$ whose cardinality tends to infinity with $x$. This subset is constructed by imposing conditions on the primes $p\leq x.$ 

By Lemma \ref{lem1} with

$$\epsilon=\min{\set{\frac{1}{2\nu f(j)} : j\in \{1,\dots,k-1\}\cup\{K\}}},$$

\vspace{2mm}
	
\noindent
$m=K,\ x_1=c_1/f(K)$ and $x_i=c_i/f(i-1)$ for $i=2,\ldots,k$, we find pairwise coprime positive integers $a_1,a_2,\ldots,a_k$ such that $(a_i,K)=1$ for any $i=1,\ldots,k$ and  
	
\begin{equation}\label{choice}
\begin{gathered}
\frac{1}{f(K)}\left(c_1-\frac{1}{2\nu}\right)<f(a_1)<\frac{1}{f(K)}\left(c_1+\frac{1}{2\nu}\right)\\
\frac{1}{f(i-1)}\left(c_i-\frac{1}{2\nu}\right)<f(a_i)<\frac{1}{f(i-1)}\left(c_i+\frac{1}{2\nu}\right)\quad \text{for} \quad i=2,\ldots,k.
\end{gathered}
\end{equation}
	
\vspace{2mm}
	
\noindent
We use these integers $a_1,a_2,\ldots,a_k$ to create the conditions that will be imposed on the primes $p\leq x.$ Particularly, in the rest of the proof we work with a subset of the primes $p\leq x$ that satisfy the linear congruences
	
\begin{equation}\label{lc}
\begin{gathered}
p+1\equiv K \mod{K^2}\\
p+i\equiv a_i \mod{a_i^2}\quad \text{for} \quad i=1,\ldots,k.
\end{gathered}
\end{equation}

\vspace{2mm}

\noindent
We shall define the precise subset of primes later in this proof. For the moment we make some observations concerning the considered prime numbers. First of all, note that such primes can be found, because the numbers $K,a_1,...,a_k$ are mutually coprime by construction, and because $(a_i-i,a_i^2)=1$ for any $i\in \set{1,\ldots,k}.$ To see the second claim, note that if there were an index $j\in \set{1,\ldots,k}$ such that $(a_j-j,a_j^2)>1,$ then there would also exist a prime $r$ that divides both $a_j-j$ and $a_j$. Then $r$ would divide $a_j-(a_j-j)=j$ and $a_j$, but this contradicts the fact that $(a_j,K)=1.$

\begin{claim}\label{cl1}
If $p$ satisfies the linear congruences (\ref{lc}), then $a_1K$ divides $p+1$, $(i-1)a_i$ divides $p+i$ for $i=2,\ldots,k$, and we also have that $((p+1)/(a_1K),a_1)=((p+1)/K,K)=1$ and that $((p+i)/((i-1)a_i),a_i)=((p+1)/(i-1),i-1)=1$ for $i=2,\ldots,k$. 
\end{claim}

\begin{proof}[Proof of Claim 1.]
Since $(a_1,K)=1$, $K|(p+1)$ and $a_1|(p+1)$, it directly follows that $a_1K$ divides $p+1$. The first linear congruence in (\ref{lc}) also implies that $((p+1)/K,K)=1.$ We continue by writing

\begin{equation}\label{help}
\frac{p+i}{i-1}=\frac{p+1}{(i-1)^2}\cdot(i-1)+1
\end{equation} 

\noindent
for $i=2,\ldots,k$. Since $K\mid (p+1),$ it is true that $(i-1)^2\mid (p+1)$ and so (\ref{help}) yields that $(i-1)|(p+i)$ and that $((p+i)/(i-1),i-1)=1.$ Moreover, when $i=2,\ldots,k$, we have that $a_i|(p+i)$ and that $(a_i,K)=1$, where the latter implies that $(a_i,i-1)=1$. Therefore, $a_i(i-1)$ divides $p+i$.

Finally we observe that the linear congruences of (\ref{lc}) imply that $((p+i)/a_i,a_i)=1$ for $i=1,\ldots,k.$ Thus, $((p+1)/(a_1K),a_1)=((p+i)/(a_i(i-1)),a_i)=1$ for $i=2,\ldots,k$ and the proof of the claim is finished. 
\end{proof}

Now we define the subset of primes that we will be working with. Let $M:=(Ka_1\cdots a_k)^2$. Because of the Chinese Remainder Theorem the primes $p\leq x$ that satisfy the linear congruences (\ref{lc}) are precisely those primes lying in some reduced residue class $N\mod{M}$ for a positive integer $N$ which is coprime to $M.$ For simplicity we set

$$\delta(p):=\frac{p+1}{a_1K}\prod_{i=2}^k\left(\frac{p+i}{a_i(i-1)}\right)$$

\noindent
and define the set of primes

$$\CS(x):=\set{p\leq x:\ p\equiv N\mod{M},\ P^-(\delta(p))>x^{\alpha},\ {\mu}^2(\delta(p))=1},$$

\vspace{2mm}

\noindent
where $\mu$ denotes the M\"{o}bius function and $\alpha$ is a fixed real number of $(0,1)$ such that

\begin{equation}\label{sie}
\#\set{p\leq x:\ p\equiv N\mod{M},\ P^-(\delta(p))>x^{\alpha}}\gg\frac{x}{(\log x)^{k+1}}.
\end{equation} 

\noindent
The existence of that real number $\alpha$ and the validity of (\ref{sie}) are proved in the next section. For now we accept this inequality whose implicit constant depends on $\alpha,k$ and $a_1,a_2,\ldots,a_k.$

\begin{claim}
If $p\in \CS(x)$, then, as $x\rightarrow +\infty$, we have that

\begin{eqnarray*}
f\left(\frac{p+1}{a_1K}\right)\!\!\!&=&\!\!\!1+o(1) \quad \text{and}\\
f\left(\frac{p+i}{a_i(i-1)}\right)\!\!\!&=&\!\!\!1+o(1) \quad \text{for} \quad i=2,\ldots,k.
\end{eqnarray*}

\end{claim}

\begin{proof}[Proof of Claim 2.]
Let $p\in \CS(x)$. When $x$ is large enough we have that

\begin{equation}\label{final}
f\left(\frac{p+i}{a_i(i-1)}\right)=\prod_{\substack{q\ \text{prime,}\ q>x^{\alpha}\\q\mid \frac{p+i}{a_i(i-1)}}}f\left(q\right)=\exp\Bigg\{O\Bigg(\sum_{\substack{q\ \text{prime,}\ q>x^{\alpha}\\q\mid \frac{p+i}{a_i(i-1)}}}|f(q)-1|\Bigg)\Bigg\}.
\end{equation}

\vspace{2mm}

\noindent
The last equality was obtained by an application of Taylor's theorem. The applicability of Taylor's theorem can be justified by the fact that $f(p)\rightarrow 1$ as $p\rightarrow \infty$, which implies that $|f(q)-1|\leq 1/2$ for the primes $q>x^{\alpha}$ with $x$ large enough. For any $i\in\set{2,\ldots,k}$ the number of distinct prime factors of the $x^{\alpha}$-rough part of $(p+i)/(a_i(i-1))$ is at most

$$\frac{\log \left(\frac{p+i}{a_i(i-1)}\right)}{\log (x^{\alpha})}\leq \frac{\log (x+i)}{\alpha \log x}\ {\ll}_{\alpha}\ 1.$$

\vspace{2mm}

\noindent
This trivial estimate implies that

$$\sum_{\substack{q\ \text{prime,}\ q>x^{\alpha}\\q\mid \frac{p+i}{a_i(i-1)}}}|f(q)-1|\ {\ll}_{\alpha}\ \sup{\set{|f(q)-1|,\ q\ \text{a prime with}\ q>x^{\alpha}}}=o(1)$$

\noindent		
as $x\rightarrow+\infty$. Thus, relation (\ref{final}) gives $f((p+i)/(a_i(i-1)))=\exp(o(1))=1+o(1)$ as $x\rightarrow+\infty.$ It is clear that the argument above leads to the same conclusion for the values $f((p+1)/(a_1K))$.
\end{proof}

We can now combine Claims 1 and 2 to show that for $i\in\set{1,\ldots,k}$, if $p\in\CS(x),$ then $f(p+i)$ is close to $c_i$. From Claim 1 and the multiplicativity of $f$ it follows that 
	
\begin{equation*}
\begin{gathered}
f(p+1)=f(K)f(a_1)f\left(\frac{p+1}{a_1K}\right)\\
f(p+i)=f(i-1)f(a_i)f\left(\frac{p+i}{a_i(i-1)}\right)\quad \text{for} \quad i=2,\ldots,k.
\end{gathered}
\end{equation*}

\noindent
Making use of (\ref{choice}), we obtain that
	
\begin{equation*}
\begin{gathered}
\left(c_1-\frac{1}{2\nu}\right)f\left(\frac{p+1}{a_1K}\right)<f(p+1)<\left(c_1+\frac{1}{2\nu}\right)f\left(\frac{p+1}{a_1K}\right)\\
\left(c_i-\frac{1}{2\nu}\right)f\left(\frac{p+i}{a_i(i-1)}\right)<f(p+i)<\left(c_i+\frac{1}{2\nu}\right)f\left(\frac{p+i}{a_i(i-1)}\right)\quad \text{for} \quad i=2,\ldots,k.
\end{gathered}
\end{equation*}
	
\vspace{2mm}

\noindent
So, when $x$ is large enough, these inequalities and Claim 2 yield that

\begin{equation*}
c_i-1/\nu<f(p+i)<c_i+1/\nu\quad \text{for} \quad i=1,\ldots,k.
\end{equation*} 

\vspace{2mm}

This means that $\CS(x)\subset S(x,\nu)$ for $x$ sufficiently large and so the proof is reduced to showing that $\#\CS(x)\rightarrow\infty$ as $x\rightarrow\infty$. To this end, we show that

\begin{equation}\label{lastin}
\#\CS(x)\gg x(\log x)^{-(k+1)}
\end{equation}

\noindent
and once this estimate is proven, the proof of the theorem will be complete.

We begin the proof of (\ref{lastin}) by taking $x$ with $x^{\alpha}\geq k-1$. Then a prime $q>x^a$ can divide at most one of the $k$ fractions of $\delta(p).$ Indeed, in the opposite case there would exist distinct indices $i,j\in \set{1,\ldots,k}$ such that $q\mid (p+i)$ and $q\mid (p+j).$ In this case $q\mid|i-j|$ and so, $x^{\alpha}<q\leq |i-j|\leq k-1$, which is a contradiction. Therefore, if $q^2\mid \delta(p)$ for a prime $q>x^{\alpha}$, there exists an $i\in \set{1,\ldots,k}$ such that $q^2\mid (p+i).$ This means that the number of primes $p\leq x$, for which $p\equiv N\mod{M}$ and $q^2\mid \delta(p)$ for some prime $q>x^{\alpha}$, is at most

\begin{equation}\label{count}
\begin{split}
\sum_{i=1}^k&\#\set{p\leq x: p\equiv N\mod{M},\ q^2\mid (p+i)\ \text{for some prime}\ q>x^{\alpha}}\leq\\
\sum_{i=1}^k&\#\set{n\leq x+i: q^2\mid n\ \text{for some prime}\ q>x^{\alpha}}\leq (x+k)\sum_{q>x^{\alpha}}\frac{1}{q^2}\ {\ll}_{\alpha,k}\ x^{1-\alpha}.
\end{split} 
\end{equation}

By (\ref{sie}) and (\ref{count}) we deduce that the number of the primes of $\CS(x)$ is $\gg x(\log x)^{-(k+1)}-x^{1-\alpha}\gg_{\alpha} x(\log x)^{-(k+1)}$ for $x$ large enough in terms of $\alpha$. In particular, we see that $\#\CS(x)\to\infty$ when $x\to\infty$, thus completing the proof of Theorem \ref{prop}.
\end{proof}

\section{Proof of Inequality (\ref{sie})}

\vspace{1mm}

To complete the proof of Theorem \ref{prop}, it remains to prove inequality (\ref{sie}). Its proof will be based upon the Fundamental Lemma of Sieve Methods.

\vspace{1mm}

\begin{thm}[The Fundamental Lemma of Sieve Theory]\label{FL}
Let $\CA$ be a finite set of integers and let $\CP$ be a set of primes. We define

\begin{eqnarray*}
A_d\!\!\!&:=&\!\!\!\#\set{a\in\CA:a\equiv 0\mod d},\\
P(y)\!\!\!&:=&\!\!\!\prod_{p\leq y,\ p\in\CP}p,\\
S(\CA,\CP,y)\!\!\!&:=&\!\!\!\#\set{a\in\CA:(a,P(y))=1}.
\end{eqnarray*}

\vspace{2mm}

\noindent
If there exists a non-negative multiplicative function $v$, some real number $X$, remainder terms $r_d$ and positive constants $\kappa$ and $C$ such that $v(p)<p$ for all $p\mid P(y)$,

\begin{eqnarray*}
A_d\!\!\!&=&\!\!\!X\cdot\frac{v(d)}{d}+r_d \quad \text{for all} \quad d\mid P(y) \quad \text{and}\\
\prod_{\substack{p\in \CP,\\w_1<p\leq w_2}}\left(1-\frac{v(p)}{p}\right)^{-1}\!\!\!&\leq&\!\!\! \left(\frac{\log w_2}{\log w_1}\right)^{\kappa}\left(1+\frac{C}{\log w_1}\right) \quad \text{for} \quad 2\leq w_1\leq w_2, 
\end{eqnarray*}

\noindent
then, for $\CA,X,y$ and $u\geq 1$ we have that

$$S(\CA,\CP,y)=X\prod_{p\leq y, p\in\CP}\left(1-\frac{v(p)}{p}\right)(1+O_{\kappa, C}(u^{-u/2}))+O\Bigg(\sum_{\substack{d\mid P(y),\\d\leq y^u}}|r_d|\Bigg).$$
\end{thm}

\begin{proof}
For a proof of this theorem the reader is advised to see the proof of Theorem $18.11$ in \cite{dimb}.
\end{proof}

The notation, which is used below, is the same as in the proof of Theorem \ref{prop}.

\vspace{1mm}

\begin{proof}[Proof of (\ref{sie})]
Like above, $p$ is a prime such that $p\leq x$ and $p\equiv N\mod{M}$. We prove that $(\delta(p),M)=1.$ This will mean that the integers, that we will be sieving, are already coprime to $M$. Thus, it suffices to sieve out the rest of the primes, so we can assume that $d$ is coprime to $M$.

Indeed, let $q$ be a prime dividing $\delta(p).$ It follows that $q\mid (p+j)/(a_j(j-1))$ for some $j\in\set{2,\ldots,k}$ or that $q\mid (p+1)/(a_1K).$ 

If $q\leq k+1,$ then $q\mid K$. In this case we can only have $q\mid (p+j)/(a_j(j-1))$ for some $j\in\set{2,\ldots,k}$, because we have already proved that $(p+1)/K$ and $K$ are coprime when $p\equiv N\mod M$, that is when $p$ satisfies the congruences of (\ref{lc}) (see Claim \ref{cl1} above). Now observe that $q$ divides $p+1$, since $q\mid K$ and $p+1\equiv K\mod{K^2}$. Thus, when $q\leq k+1$ we have that $q\mid (p+j)-(p+1)=j-1$. However, it was shown that $(p+j)/(j-1)$ and $j-1$ are coprime when $p\equiv N\mod M$. Consequently, $q$ cannot divide both $(p+j)/(a_j(j-1))$ and $j-1$. Hence, $q>k+1.$ 

If $q$ divides $a_{\ell}$ for some $\ell\in \set{1,\ldots,k}$, then $q\mid p+\ell,$ since $p+\ell\equiv a_{\ell}\mod{{a_{\ell}}^2}.$ Moreover, $q\mid(p+j)/a_j$ for some $j$, because $q$ divides $\delta(p).$ Since $(p+j)/a_j$ and $a_j$ are coprime and we know that $q\mid a_{\ell}$, we must have that $j\neq \ell.$ Now, we have $q|p+j$ and $q|p+\ell$, whence $q|(j-\ell).$ In particular, $q\le k-1$, but this is impossible since $q|a_{\ell}$. Thus, $(q,a_j)=1$ for all $j\in \set{1,\ldots,k}$. This proves our claim that $(\delta(p),M)=1$ when $p\equiv N\mod M.$ 

Now we can rewrite the left-hand side of (\ref{sie}) as $S(\CA,\CP,y)$ with $\CA=\{\delta(p):p\leq x,\ p\equiv N\mod M\}$, $\CP=\{p\leq y\}$ and $y=x^{1/u}$. We shall verify the hypotheses of Theorem \ref{FL}. Let $d\mid \delta(p)$. We then have that $(d,M)=1.$ Hence, a prime $p\leq x$ is counted by $A_d$ if and only if $p\equiv N \mod M$ and $p\equiv a \mod d$, where $a$ is such that $(a,d)=1$ and $(a+1)\cdots(a+k)\equiv 0\mod d$. Let us write $g(d)$ for the total number of such $a$'s mod $d$. For each fixed $a$, the number of admissible primes is $\pi(x)/\phi(Md) + O(E(Md))$, where

$$E(q)=\max_{\substack{1\leq b<q\\(b,q)=1}}\abs{\pi(x;q,b)-\frac{\pi(x)}{\phi(q)}}.$$

\noindent
Consequently, 

$$A_d=\frac{g(d)}{\phi(d)}\cdot\frac{\pi(x)}{\phi(M)} + O\Big(g(d)E(Md)\Big).$$

\noindent
Note that $g(p)=k$ for each prime $p > k+1$. Hence, $g(d)=k^{\omega(d)}.$ So, Theorem \ref{FL} implies that 

\begin{equation*}
\begin{split}
\#\set{p\leq x: p\equiv N\mod{M},\ P^-(\delta(p))>y}&=\frac{\pi(x)}{\phi(M)}\prod_{P^+(M)<p\leq y}\left(1-\frac{k}{p-1}\right)(1+O(u^{-u/2}))\\
&+O\Bigg(\sum_{d\leq y^u}k^{\omega{(d)}}E(Md)\Bigg)
\end{split}
\end{equation*}

\vspace{2mm}

\noindent
and this holds for $x,y,u=:1/3\alpha\geq 1.$ For $\alpha$ small enough, the main term is ${\gg}_{k,M}\frac{x}{\log x(\log y)^k}.$ 

We focus on bounding the sum in the error term. The Bombieri-Vinogradov theorem guarantees the existence of a positive constant $B=B(k)$ such that

$$\sum_{q\,<\frac{x^{1/2}}{(\log x)^B}}\!\!\!\!E(q)\ll \frac{x}{(\log x)^{2(k^2+k+2)}}.$$

\noindent
Therefore, if we choose $y=x^{\alpha}$, then for $x$ large enough we have that $My^{\frac{1}{3\alpha}}<x^{1/2}/(\log x)^B$ and so, by the trivial fact $k^{\omega(q)}<k^{\omega(Mq)}$ and the Cauchy-Schwarz inequality we obtain that the sum in the big-O term above is bounded by

\begin{eqnarray*}
	\sum_{q\,<\frac{x^{1/2}}{(\log x)^B}}\!\!\!\!k^{\omega{(q)}}E(q)\!\!\!&\leq&\!\!\!\Bigg(\sum_{q\,<\frac{x^{1/2}}{(\log x)^B}}\!\!\!\!E(q)\Bigg)^{1/2}
	\Bigg(\sum_{q<x^{1/2}}k^{2\omega{(q)}}E(q)\Bigg)^{1/2}\\&\ll&\!\!\!\frac{x^{1/2}}{(\log x)^{k^2+k+2}}\Bigg(\sum_{q<x^{1/2}}k^{2\omega{(q)}}E(q)\Bigg)^{1/2}.
\end{eqnarray*}

\noindent
The trivial estimates $\pi(x;q,a)\ll x/q,\ \pi(x)\!\leq \!x$ and $\phi(q)\gg q/(\log \log q)\gg q/(\log \log x)$ imply that

$$\left|\pi(x;q,a)-\frac{\pi(x)}{\phi(q)}\right|\leq \max{\set{\pi(x;q,a),\frac{\pi(x)}{\phi(q)}}}\ll \frac{x\log \log x}{q}.$$

\noindent
Consequently,

\begin{eqnarray*}
	\frac{1}{x\log\log x}\sum_{q<x^{1/2}}k^{2\omega{(q)}}E(q)\!\!\!&\ll&\!\!\!\sum_{q\leq x}\frac{k^{2\omega(q)}}{q} \leq \prod_{p\leq x}\left(1+\frac{k^2}{p}+\frac{k^2}{p^2}+\ldots\right)\\
	&=&\!\!\!\prod_{p\leq x}\left(1+\frac{k^2}{p-1}\right)
	{\ll}_k\,(\log x)^{k^2}
\end{eqnarray*}

\noindent
and so,

$$\sum_{q\,<\frac{x^{1/2}}{(\log x)^B}}\!\!\!\!k^{\omega{(q)}}E(q)\ {\ll}_k\,\frac{x}{(\log x)^{k+2}}.$$

\noindent
This means that the error term above is ${\ll}_k\ x/(\log x)^{k+2}.$ 

With the choice $u=1/3\alpha$ for $\alpha$ small enough in terms of $k$, we have that the main term is ${\gg}_{u,k,M}\ x/(\log x)^{k+1}.$ Therefore,

$$\#\set{p\leq x: p\equiv N\mod{M},\ P^-(\delta(p))>x^{\alpha}}{\gg}_{\alpha,k,M}\,\frac{x}{(\log x)^{k+1}}$$

\vspace{2mm}

\noindent
for $x$ large enough in terms of $k$. So, the proof of (\ref{sie}) is complete.
\end{proof}

\vspace{1mm}


\bibliographystyle{alpha}

\end{document}